\documentclass [12pt,a4paper]{amsart}
\usepackage{stmaryrd}
\usepackage{textcomp}
\usepackage{amsmath}
\usepackage{cases}
\usepackage{amscd}
\usepackage{ifthen}

\usepackage{epsfig}
\usepackage{graphicx}
\usepackage{verbatim}
\usepackage{amssymb}
\usepackage{amsfonts}
\usepackage{amsbsy}
\usepackage{amsthm}
\usepackage{amstext}
\usepackage{amsopn}
\usepackage[dvips]{color}
\usepackage{indentfirst}
\usepackage{mathrsfs}
\usepackage{enumerate}

\newtheorem{thm}{Theorem}[section]

\newtheorem{cor}[thm]{Corollary}
\newtheorem{lem}[thm]{Lemma}
\newtheorem{prop}[thm]{Proposition}

\newtheorem*{thmA}{Theorem A}

\newenvironment{rem}{%
\bigskip
\noindent
\textsl{{\sl Remark. }}}{\bigskip}

\newenvironment{pf}[1][]{%
 \vskip 3mm
 \noindent
 \ifthenelse{\equal{#1}{}}%
  {{\slshape Proof. }}%
  {{\slshape #1.} }%
 }%
{\qed\bigskip}

\newcommand{\A}{{\mathcal A}}

\renewcommand{\arg}{{\operatorname{arg}\,}}
\newcommand{\arctanh}{{\operatorname{arctanh}\,}}

\newcommand{\C}{{\mathbb C}}

\newcommand{\D}{{\mathbb D}}

\newcommand{\bD}{{\overline{\mathbb D}}}

\newcommand{\F}{{\mathcal F}}

\renewcommand{\Im}{{\operatorname{Im}\,}}

\newcommand{\K}{{\mathcal K}}

\newcommand{\Li}{{\operatorname{Li}}}

\newcommand{\Cara}{{\mathcal P}}

\newcommand{\R}{{\mathbb R}}

\renewcommand{\Re}{{\operatorname{Re}\,}}
\newcommand{\es}{{\mathcal S}}

\renewcommand{\SS}[1]{{\mathcal{S}^*_{#1}}}

\newcommand{\inv}{^{-1}}

\usepackage{color}

\newcounter{minutes}
\divide\time by 60
\newcounter{hours}
\multiply\time by 60
\addtocounter{minutes}{-\time}
\begin{document}
\bibliographystyle{amsplain}
\title[Open door functions]%
{Analytic and geometric properties of \\
open door functions}

\author[M. Li]{Ming Li}
\address{Graduate School of Information Sciences \\
Tohoku University \\
Aoba-ku, Sendai 980-8579, Japan}
\email{li@ims.is.tohoku.ac.jp}
\author[T. Sugawa]{Toshiyuki Sugawa}
\address{Graduate School of Information Sciences \\
Tohoku University \\
Aoba-ku, Sendai 980-8579, Japan}
\email{sugawa@math.is.tohoku.ac.jp}

\begin{abstract}
In this paper, we study analytic and geometric properties of
the solution $q(z)$ to the differential
equation $q(z)+zq'(z)/q(z)=h(z)$ with the initial condition $q(0)=1$
for a given analytic function $h(z)$ on the unit disk $|z|<1$ in the
complex plane with $h(0)=1.$
In particular, we investigate the possible largest constant $c>0$
such that the condition $|\Im[zf''(z)/f'(z)]|<c$ on $|z|<1$ implies
starlikeness of an analytic function $f(z)$ on $|z|<1$ with
$f(0)=f'(0)-1=0.$
\end{abstract}
\maketitle

\section{Introduction}
We denote by $\A$ the class of holomorphic functions on the unit disk
$\D=\{z: |z|<1\}$ of the complex plane $\C.$
Let $\A_0$ denote the subclass of $\A$ consisting of functions
$p$ with $p(0)=1.$
Let $\A_1$ be the class of functions of the form $zp(z)$ for $p\in\A_0.$
In other words, $f\in\A_1$ if and only if $f\in\A$ and $f(0)=f'(0)-1=0.$
We say that a function $f\in\A$ is {\it subordinate} to another $g\in\A$
and write $f\prec g$ or $f(z)\prec g(z)$
if $f=g\circ\omega$ for a function $\omega\in\A$ such that
$\omega(0)=0$ and $|\omega|<1.$
When $g$ is univalent, $f\prec g$ precisely when $f(0)=g(0)$
and $f(\D)\subset g(\D).$

The set of functions $q\in\A_0$ with $\Re q>0$ is called the
Carath\'eodory class and will be denoted by $\Cara.$
It is well recognized that the function $q_*(z)=(1+z)/(1-z)$
(or its rotation $q_*(e^{i\theta}z)$) maps the unit disk univalently onto
the right half-plane and is extremal in many problems.
A function $f\in\A_1$ is called {\it starlike} if $f$ maps $\D$ univalently
onto a starlike domain with respect to the origin.
Likewise, a function $f\in\A_1$ is called {\it convex} if $f$ maps $\D$
univalently onto a convex domain.
We denote by $\es^*$ and $\K$ the classes of starlike and convex functions,
respectively.
It is well known that $f\in\A_1$ is starlike precisely when
$q(z)=\psi_f(z):=zf'(z)/f(z)$ belongs to $\Cara$ and that
$f\in\A_1$ is convex precisely when $h(z)=\varphi_f(z):=1+zf''(z)/f'(z)$
belongs to $\Cara$ (see, for instance, \cite{Duren:univ}).
Note here the relation $h(z)=q(z)+zq'(z)/q(z).$
We also note that $f(z)$ is convex if and only if $zf'(z)$ is starlike
for $f\in\A_1.$
For a given $h\in\A_0,$ we always find a function $f\in\A_1$
with $1+zf''/f'=h.$
Indeed, by integrating the relation $(\log f')'=f''/f'=(h-1)/z,$
we obtain
\begin{equation}\label{eq:convex}
f'(z)=\exp\left( \int_0^z\frac{h(t)-1}t dt \right).
\end{equation}
By integrating the above, we get the desired $f\in\A_1.$
Similarly, replacing $f'(z)$ by $f(z)/z$ in \eqref{eq:convex}, we obtain
a representation of $f\in\A_1$ satisfying $zf'/f=h.$

It is obvious that a convex function is starlike.
This means analytically that $h=q+zq'/q\in\Cara$ implies
$q\in\Cara$ for $q\in\A_0.$
In other words, $q+zq'/q\prec q_*$ implies $q\prec q_*.$

One can observe that the function
$$
h_*(z):=q_*(z)+\frac{zq_*'(z)}{q_*(z)}
=\frac{1+z}{1-z}+\frac{2z}{1-z^2}
=\frac{1+4z+z^2}{1-z^2}
$$
maps the unit disk onto the complex plane $\C$ slit along
the two half-lines $\pm iy,~y\ge\sqrt3.$
The following was proved by Mocanu \cite{Mocanu86} and later
extended by Miller and Mocanu \cite{MM97BB}
(see also \cite{MM:ds}).

\begin{thmA}[Open Door Lemma]\label{Thm:ODL}
Suppose that a function $q\in\A_0$ satisfies the subordination condition
$$
q(z)+\frac{zq'(z)}{q(z)}\prec h_*(z)=q_*(z)+\frac{zq_*'(z)}{q_*(z)}.
$$
Then $q(z)\prec q_*(z).$
\end{thmA}

In particular, if a function $f\in\A_1$ satisfies the
subordination $1+zf''/f'\prec h_*,$ then $f$ is starlike.
Since the slit domain $h_*(\D)$ contains the parallel strip
$|\Im w|<\sqrt3,$ we obtain the following result as a corollary.

\begin{cor}\label{cor:sqrt3}
If a function $f\in\A_1$ satisfies the condition
$$
\left|\Im\left[\frac{zf''(z)}{f'(z)}\right]\right|<\sqrt3,
\quad z\in\D,
$$
then $f$ is starlike.
\end{cor}

%Next we give an example to show that this slit domain is not the best one.
%\begin{eg}
%For $q(z)=\frac{(1+\sqrt{5}i)z}{1-e^{-(1+\sqrt{5}i)z}}$ and  $h(z)=1+(1+\sqrt{5}i)z$.
%Obviously $q,h\in\mathcal{H}[1,1]$ satisfy $q(0)=1$,
%$\frac{zq'(z)}{q(z)}+q(z)=h(z)$, and $h(\mathbb{U})\not\subset V(-\sqrt{3},\sqrt{3})$, $\Re q(z)>0$ for $z\in\mathbb{U}$.
%\end{eg}

We recall a notion of strong starlikeness.
A function $f\in\A_1$ is called {\it strongly starlike of order $\alpha$}
for an $0<\alpha$ if $|\arg [zf'(z)/f(z)]|<\pi\alpha/2$ for $z\in\D.$
We denote by $\SS{\alpha}$ the class of strongly starlike functions in $\A_1$
of order $\alpha.$
Obviously, we have $\SS{1}=\es^*.$
For geometric characterizations of strongly starlike functions, see
\cite{SugawaDual} and references therein.

In the present paper, we try to find or estimate the best possible constant
$\gamma>0$ such that the condition $|\Im [zf''(z)/f'(z)]|<\gamma$
implies $f\in\SS{\alpha}.$
%For $h\in\A_0,$ we denote by $\delta(h)$ the largest $\gamma>0$ such that
%the parallel strip $|\Im w|<\gamma$ is contained in $h(\D).$
%If there is no such a parallel strip, we put $\delta(h)=0.$
%In other words,
More precisely, the number is defined as $\gamma(\SS{\alpha}),$ where
$$
\gamma(\F)=\sup\big\{\gamma\ge0: \varphi_f(\D)\subset W_\gamma
~\text{implies}~ f\in\F~\text{for}~ f\in\A_1\big\}
$$
for a subset $\F$ of $\A_1$ and
$$
W_\gamma=\{w\in\C: |\Im w|<\gamma\}
$$
is a parallel strip of width $2\gamma.$
We recall that $\varphi_f=1+zf''/f'.$

We will show the following estimates of $\gamma(\SS{\alpha}).$
See also Figure \ref{fig:gamma} below for the graphs of 
$\gamma(\SS{\alpha})$ the upper and lower bounds.

\begin{thm}\label{thm:ss}
Let $0<\alpha<1.$
Then
$$
\sqrt3\alpha<
\frac{\alpha+(1+\alpha)\sin(\pi\alpha/2)}{\sqrt{1+2\sin(\pi\alpha/2)}}
<\gamma(\SS{\alpha})<\frac{\sqrt3\pi\alpha}{\sqrt{3+\alpha}}.
$$
\end{thm}

We remark that a similar (but not better) result can be found
at \cite[Theorem 1.6]{KPS04}.
Mocanu \cite[Corollary 1.1]{Mocanu86A} showed that $\gamma(\SS{\pi/4})\ge 1.$
Our estimate gives $\gamma(\SS{\pi/4})>(2+3\sqrt2)/4\sqrt{1+\sqrt2}
=1.0044319\dots.$
Note that the lower bound in this theorem tends to $\sqrt3$ as $\alpha\to1,$
which agrees with Corollary \ref{cor:sqrt3}.
When $\alpha=1,$ we can slightly improve the upper bound in
the last theorem.

\begin{thm}\label{thm:starlike}
$\sqrt3\le \gamma(\es^*)<2.5.$
%\dfrac{\sqrt3\pi}{2}.$
\end{thm}

Though it is difficult to compute the exact value of $\gamma(\es^*),$
the next result gives us a way to compute it numerically.

\begin{thm}\label{thm:char}
Let $\theta_c=2\arctan(e^{2/c})\in(\pi/2,\pi)$ for $c>0.$
Let $F(c)=v(1)$ for $c\ge0,$ where $v(t)$ is the solution to the
initial value problem of ordinary differential equation
$$
v(t)+\frac{tv'(t)}{v(t)}=1+\frac c2\log\frac{1+te^{i\theta_c}}{1-te^{i\theta_c}},
\quad v(0)=1.
$$
Then $\gamma(\es^*)=\pi c_0/4,$ where $c_0$ is the smallest positive number
such that $\Re F(c_0)=0.$
\end{thm}

We remark that $F(c)=q_c(e^{i\theta_c}),$ where $q_c$ is given in Section 3.

The organization of the present paper is as follows.
In Section 2, we investigate geometric properties of the solution
$q$ to the differential equation $q+zq'/q=h$ for a given $h\in\A_0.$
We believe that our observation will be helpful for other kinds of
problems concerning the subordination of $\psi_f$ and $\varphi_f.$
In order to estimate the quantity $\gamma(\SS{\alpha}),$ 
we study in Section 3 the extremal case when 
$h=h_c=1+c\cdot\arctanh,$ 
which maps $\D$ onto the parallel strip $W_{\pi c/4}$ for $c>0.$
We will show that the solution $q_c$ to the differential equation
$q_c+zq_c'/q_c=h_c$ maps $\D$ univalently onto a smooth Jordan
domain if $c$ is not very large (Theorem \ref{thm:qc}).
Lemma \ref{lem:ca} will describe $\gamma(\SS{\alpha})$
in terms of the above solutions $q_c.$
Section 4 is devoted to the proof of the main theorems.
The last section gives concluding remarks on numerical experiments.
By using Mathematica, we can generate a graph of the function
$\alpha\mapsto\gamma(\SS{\alpha})$ and some approximation
of the value $\gamma(\es^*)$ though there is no rigorous error
estimate for the present experiments.

\section{Some observations of Open Door Lemma}

Throughout the present paper, for simplicity,
we consider only functions in $\A_0$ for the Open Door Lemma.
For the most general version of Open Door Lemma, the reader should
consult the monograph \cite{MM:ds} by Miller and Mocanu.

The following result is contained in
Theorems 3.2i and 3.4b of \cite[p.~97, p.~124]{MM:ds}.

\begin{lem}\label{lem:MM}
Let $h\in\A_0$ map $\D$ univalently onto a convex domain.
Suppose that the differential equation
\begin{equation}\label{eq:qh}
q(z)+\frac{zq'(z)}{q(z)}=h(z)
\end{equation}
has an analytic solution $q$ with $\Re q>0.$
Then $q$ is univalent and $q\prec h.$
If $p\in\A_0$ satisfies
\begin{equation}\label{eq:ph}
p(z)+\frac{zp'(z)}{p(z)}\prec h(z),
\end{equation}
then $p\prec q$ and $q$ is the best dominant.
\end{lem}

The lemma immediately yields the following corollary.

\begin{cor}
Let $h\in\A_0$ be a univalent function with convex image containing
the parallel strip $W_\gamma.$
If \eqref{eq:qh} has an analytic solution $q$ with $\Re q>0,$
then $\gamma(\es^*)\ge\gamma.$
\end{cor}

It is, in general, not easy to analyse the solution to the
differential equation \eqref{eq:qh} for a given $h.$
Therefore, practically, we start from a function $q$ with $\Re q>0$
and look at the image of $\D$ under the function $h$ defined by
\eqref{eq:qh}.
If it is a convex domain containing $W_\gamma,$ then the last corollary
implies $\gamma(\es^*)\ge\gamma.$
Therefore, to make a suitable choice of $q,$ it is helpful to observe
the boundary behaviour of the solution to the equation \eqref{eq:qh}
for the targeted $h.$

For $q\in\A_0$ and $\theta\in\R,$ we define
%$\beta_+(\theta)$ and
$\beta(\theta)=\beta(\theta;q)$ (modulo $2\pi$ rigorously speaking) by
$$
\beta(\theta)=\lim_{z\to e^{i\theta}}\left[\arg q'(z)+\theta+\frac\pi2\right]
$$
if it exists.
When $q$ and $q'$ extend continuously to $\{z\in\bD: |z-e^{i\theta}|<\delta\}$
for some $\delta>0$ and $q'(e^{i\theta})\ne0,$ one has
$$
\lim_{t\to\theta}\arg\frac{d}{dt}q(e^{it})
=\lim_{t\to\theta}\left[\arg q'(e^{it})+t+\frac\pi2\right]
=\beta(\theta).
$$
Thus, $\beta(\theta)$ means the argument of a tangent vector
of the boundary curve $q(e^{it})$ at $t=\theta.$
Even if the above limit does not exist, the following limits may exist:
$$
\beta_\pm(\theta)
=\lim_{t\to\theta\pm}\arg\frac{d}{dt}q(e^{it})
=\lim_{t\to\theta\pm}\left[\arg q'(e^{it})+t+\frac\pi2\right]
$$
for each signature.
%If both $\beta_+(\theta)$ and $\beta_-(\theta)$ exist and are equal,
%we denote the common value by $\beta(\theta).$
%For instance, if the curve $t\mapsto q(e^{it})$ is smooth and has a non-vanishing
%derivative at $t=\theta,$ then $\beta_+(\theta)=\beta_-(\theta).$
Assume that $q$ maps $\D$ univalently onto a Jordan domain $\Omega.$
$\Omega$ is said to be {\it smooth} if $\beta$ can be chosen as a continuous
function on $\R.$
This means that $\Omega$ has continuously varying tangent at each boundary
point.
For further properties of $\beta(\theta),$ see \cite[\S 3.2]{Pom:bound}.

We summarise properties of the solutions to \eqref{eq:qh}.

\begin{prop}\label{prop:prop}
Let $h\in\A_0$ and take $f\in\A_1$ so that \eqref{eq:convex} is fulfilled.
Then a unique meromorphic solution $q$ to the differential equation \eqref{eq:qh}
with $q(0)\ne0$ is given by $q(z)=zf'(z)/f(z).$
If the solution $q$ has a pole at $z=z_0,$ then its order is $1$ and
its residue is $z_0.$

%\item
Suppose that $h(z)=O(|z-\zeta|^{-\alpha})$ as 
$z$ tends to a boundary point $\zeta=e^{i\theta}\in\partial\D$ in
$\D$ for a constant $\alpha<1.$
Then the limits $\lim_{z\to\zeta}f'(z)=:f'(\zeta)$ and
$\lim_{z\to\zeta}f(z)=:f(\zeta)$ exist and $f'(\zeta)\ne0.$
If, in addition, $f(\zeta)\ne0,$ then $\lim_{z\to \zeta}q(z)=:q(\zeta)$ exists
and $q(\zeta)\ne0.$
Moreover the following hold.
\begin{enumerate}
\renewcommand{\labelenumi}{$({\roman{enumi}})$}
\item
Suppose that the finite limit $\lim_{z\to\zeta}h(z)=:h(\zeta)$ exists
and $f(\zeta)\ne0.$
Then,
$q(z)=q(\zeta)[1+\bar\zeta(h(\zeta)-q(\zeta)+o(1))(z-\zeta)]$ as $z\to\zeta.$
Moreover, if $q(\zeta)\ne h(\zeta),$ then
$$
\beta(\theta)=\arg q(\zeta)+\arg[h(\zeta)-q(\zeta)]+\frac\pi2.
$$
%Suppose that $h(z)=A+(B+o(1))(z-\zeta)^\alpha$ as $z\to\zeta$ in $\D$
%for constants $A, B, \alpha$ with $B\ne0$ and $0<\alpha<1.$
%If $f(\zeta)\ne0,$ then $q(\zeta)=\lim_{z\to \zeta}q(z)$ exists and
%$q(z)=q(\zeta)[1+(A-q(\zeta)+o(1))(z-\zeta)/\zeta]$ as $z\to\zeta.$
\item
If $h(z)=(A+o(1))(\zeta-z)^{-\alpha}$ as $z\to\zeta$ in $\D$
for constants $A\ne0$ and $0<\alpha<1,$
then $q(z)=q(\zeta)[1-\frac{\bar\zeta A}{1-\alpha}(\zeta-z)^{1-\alpha}]$
as $z\to\zeta.$
$$
\beta_\pm(\theta)=\arg q(\zeta)+\arg A-\alpha\theta
+(1\pm\alpha)\frac\pi2.
$$
\item
If $h(z)=-(A+o(1))\log(\zeta-z)$ as $z\to\zeta$ in $\D$
for a constant $A\ne0,$ then
$q(z)=q(\zeta)[1+(\bar\zeta A+o(1))(z-\zeta)\log(z-\zeta)]$ and
$$
\beta(\theta)=\arg q(\zeta)+\arg A+\frac\pi2.
$$
\end{enumerate}
\end{prop}

\begin{rem}
In the limit above, $z\to\zeta$ means that $z$ approaches $\zeta$
in $\D$ without any restriction such as radial or non-tangential limits.
\end{rem}

\begin{pf}
Observe first that any analytic function $q$ with $q(0)\ne0$
satisfying \eqref{eq:qh}
on $\D$ must have the initial value $q(0)=1$ because of $h(0)=1.$
We now show that \eqref{eq:qh} has a
meromorphic solution $q$ on $\D.$
Since $\D$ is simply connected, an analytic function
$f$ on $\D$ can be defined uniquely by \eqref{eq:convex} and the
condition $f(0)=0.$
Note that $f'(z)\ne0$ for $z\in \D.$
Therefore, $q(z)=zf'(z)/f(z)$ is a meromorphic solution
to \eqref{eq:qh} on $\D$ with at most simple poles.

It is easy to show existence of the limits when $z\to\zeta$ in $\D.$
(See the proof of (ii) to get basic ideas to do that.)
Hence, we have shown the first assertion in the proposition.

We show now assertion (i).
As for the formula of $\beta(\theta),$ one needs only to take the argument
of both sides of the identity
\begin{equation}\label{eq:arg}
%e^{i\theta}q'(e^{i\theta})=
zq'(z)
=q(z)(h(z)-q(z))
\end{equation}
and put $z=\zeta=e^{i\theta}.$

We next show assertion (ii).
Suppose that $h(z)=(A+o(1))(\zeta-z)^{-\alpha}$ as $z\to\zeta$
in $\D$ for $\alpha\in\R$ with $0<\alpha<1.$
(We will surpress the description ``in $\D$" in the rest of the proof
for brevity.)
Since $0<\alpha<1,$ we have
$$
\frac{h(z)-1}z=\left(\frac A\zeta+o(1)\right)(\zeta-z)^{-\alpha},
\quad z\to\zeta,
$$
and the limit
$$
C=\lim_{z\to\zeta}\int_0^z\frac{h(t)-1}t dt
=\int_0^\zeta\frac{h(t)-1}t dt
$$
exists.
Thus, in view of \eqref{eq:convex}, one obtains
\begin{align*}
e^{-C}f'(z)&=\exp\left[\int_\zeta^z \frac{h(t)-1}{t}dt\right] \\
&=\exp\left[\left(-K+o(1)\right)(\zeta-z)^{1-\alpha}
\right] \\
&=1+(-K+o(1))(\zeta-z)^{1-\alpha}
\end{align*}
as $z\to\zeta,$ where $K=\bar\zeta A/(1-\alpha).$
In particular, the limits $\lim_{z\to\zeta}f'(z)=e^C=:f'(\zeta),$
$\lim_{z\to\zeta}f(z)=:f(\zeta)$ exist and
$$
f(z)=f(\zeta)+f'(\zeta)
\big[(z-\zeta)+(\tfrac{K}{2-\alpha}+o(1))(\zeta-z)^{2-\alpha}\big],
\quad z\to\zeta.
$$
If $f(\zeta)\ne0,$ then $f(z)/f(\zeta)=1+q(\zeta)(z-\zeta)/\zeta
+O(|z-\zeta|^{2-\alpha})$ so that
$$
\frac{q(z)}{q(\zeta)}=1+(-K+o(1))(\zeta-z)^{1-\alpha}+\frac{1-q(\zeta)}\zeta
(z-\zeta)
$$
as $z\to\zeta.$
Thus the first part of (ii) has been shown.

To show the relation for $\beta_\pm(\theta)$ in (ii),
we note that $\arg[(\zeta-z)\inv]\to -\theta\pm\pi/2$ as $t\to\theta\pm$ for
$\zeta=e^{i\theta}.$
Since $|h(z)|\to+\infty$ as $z\to\zeta$ in this case, \eqref{eq:arg} yields
\begin{align*}
\beta_\pm(\theta)
&=\arg q(\zeta)+\lim_{t\to\theta\pm}\arg[h(e^{it})-q(e^{it})]+\frac\pi2 \\
&=\arg q(\zeta)+\lim_{t\to\theta\pm}\arg h(e^{it})+\frac\pi2 \\
&=\arg q(\zeta)+\arg A+\alpha\left(-\theta\pm\frac\pi2\right)+\frac\pi2.
\end{align*}

Finally, we show assertion (iii).
We can compute in the same way as in (ii) except for the integrals:
$$
\int_\zeta^z\log(t-\zeta)dt=(z-\zeta)[\log(z-\zeta)-1]
=(1+o(1))(z-\zeta)\log(z-\zeta)=o(1)
$$
and
$$
\int_\zeta^z(t-\zeta)\log(t-\zeta)dt=\frac{(z-\zeta)^2}2[\log(z-\zeta)-1/2]
=O\left((z-\zeta)^2\log(z-\zeta)\right)
$$
as $z\to\zeta.$
Thus the conclusion follows.
%The formula for $\beta(\theta)$ can be deduced in the same way as above.
\end{pf}

In the last proposition, the condition in (ii) means roughly that
$z=\zeta$ corresponds via the function $h(z)$ to the tip at infinity
of a sector with opening angle $\pi\alpha.$
It should also be noted that assertion (iii) can be regarded as
a limiting case of assertion (ii).

\section{An extremal case}

For $c>0,$ we define
$$
h_c(z)=1+c\,\arctanh z
=1+\frac c2\log\frac{1+z}{1-z}=1+c\sum_{n=0}^\infty\frac{z^{2n+1}}{2n+1},
\quad z\in\D,
$$
and let $q_c$ be the solution to the initial value problem of the ODE:
\begin{equation}\label{eq:q}
q_c(z)+\frac{zq_c'(z)}{q_c(z)}=h_c(z), \quad z\in\D, \quad q_c(0)=1.
\end{equation}
Note that $h_c$ maps the unit disk $\D$ onto the parallel strip
$W_{\pi c/4}=\{w: |\Im w|<\pi c/4\}.$
Let $f_c\in\A_1$ be the solution to the equation $\varphi_{f_c}=h_c.$
Namely, $f=f_c$ can be determined by \eqref{eq:convex} with $h=h_c.$
Then, $q_c(z)=zf_c'(z)/f_c(z).$
We compute
$$
\int_0^z\frac{h_c(t)-1}t dt=c \chi_2(z),
$$
where
$$
\chi_2(z)=\frac12\left[\Li_2(z)-\Li_2(-z)\right]
=\sum_{n=0}^\infty\frac{z^{2n+1}}{(2n+1)^2}
$$
is called Legendre's chi-function (see \cite[\S 1.8]{Lewin:poly})
and $\Li_2(z)$ is the dilogarithm function.
Note that
$$
\chi_2(1)=\sum_{n=0}^\infty\frac{1}{(2n+1)^2}=\frac{\pi^2}8<+\infty.
$$
Therefore, by \eqref{eq:convex}, $f_c$ is expressed by
$$
f_c(z)
=\int_0^z\exp\big[c \chi_2(t)\big]dt
=z\int_0^1\exp\big[c \chi_2(tz)\big]dt.
$$
Hence,
\begin{equation}\label{eq:1/q}
\frac1{q_c(z)}=\frac{f_c(z)}{zf_c'(z)}
=\int_0^1\exp c\big[\chi_2(tz)-\chi_2(z)\big]dt.
\end{equation}

%We recall a useful lemma due to Miller and Mocanu (see \cite[Lemma 2.2c]{MM:ds}).
%Let $E(h)=\{\zeta\in\partial\D: h(z)\to\infty ~\text{as}~ z\to\zeta
%~\text{in}~\D\}$ for an analytic function $h$ on $\D.$

%\begin{lem}\label{lem:Jack}
%Let $q$ and $h$ be non-constant analytic functions on $\D$
%with $q(0)=h(0).$
%Assume that $h$ is analytic and injective on $\bD\setminus E(h).$
%If there exist points $z_0\in\D$ and $\zeta_0\in\partial\D\setminus E(h)$
%such that $q(z_0)=h(\zeta_0)$ and that $q(\{|z|<r_0\})\subset h(\D),$
%where $r_0=|z_0|.$
%Then there exists a real number $m\ge1$ such that
%$z_0q'(z_0)=m\zeta_0h'(\zeta_0).$
%\end{lem}

We define two numbers $c_1$ and $c_*$ as the largest possible ones
with the properties
\begin{align*}
0<c<c_1\quad &\Rightarrow \quad \Re q_c>0~\text{on}~\D, \\
0<c<c_*\quad &\Rightarrow \quad q_c\prec h_c.
\end{align*}
By Corollary \ref{cor:sqrt3} and Lemma \ref{lem:MM}, we observe
that the following inequalities hold:
$$
\frac{4\sqrt3}{\pi}\le c_1\le c_*.
$$
It is also easy to show that $\Re q_{c_1}>0$ and $q_{c_*}\prec h_{c_*}.$
We are now able to show the following.

\begin{thm}\label{thm:qc}
Let $0<c\le c_*.$
Then the solution $q_c$ to \eqref{eq:q} is a non-vanishing analytic function
on the unit disk $\D$ and satisfies the relation
$\overline{q_c(z)}=q_c(\bar z)$ and the inequalities
%$q_{c}\prec h_c$ and
$$
0<q_c(-1)<|q_c(z)|<q_c(1)<+\infty,\quad z\in\D.
$$
If, in addition, $c\le c_1,$ then $q_c$ is univalent on $\D$ and
the image $q_c(\D)$ is a smooth Jordan domain in the sense that
its boundary has continuously varying tangent.
\end{thm}

\begin{pf}
By Proposition \ref{prop:prop}, we first see that $q_c$ extends meromorphically
to $\bD\setminus\{1,-1\}$ and that $q_c(z)$ has finite limits $q_c(\pm1)$
as $z\to\pm1$ in $\D.$
The symmetry property in the real axis is immediate from
uniqueness of the initial value problem for ODE.

We now prove the inequalities in the assertion.
In view of the expression \eqref{eq:1/q}, the reciprocal $1/q_c$
is analytic on $\D.$
We now look at the function $p(x)=1/q_c(x)$ for $-1<x<1.$
Since
$$
\chi_2'(x)-t\chi_2'(tx)=\sum_{n=0}^\infty\frac{1-t^{2n+1}}{2n+1}x^{2n}>0
$$
for $-1<x<1$ and $0<t<1,$ one obtains $p'(x)<0,$ which implies that
$q_c(x)$ is increasing in $-1<x<1.$

Let $u(\theta)=R(\theta)e^{i\Theta(\theta)}=q_{c}(e^{i\theta})$
with $\Theta(0)=0$ for $0\le\theta\le\pi.$
Then, we deduce from \eqref{eq:q} that
\begin{align*}
q_c(e^{i\theta})+\frac{e^{i\theta}q_c'(e^{i\theta})}{q_c(e^{i\theta})}
&=u(\theta)+\frac{u'(\theta)}{iu(\theta)}
=R(\theta)e^{i\Theta(\theta)}+\frac{R'(\theta)}{iR(\theta)}+\Theta'(\theta) \\
&=\frac c2\log{\frac{1+e^{i\theta}}{1-e^{i\theta}}}+1
=\frac c2\left(\log\cot\frac\theta2+\frac{\pi i}2\right)+1
\end{align*}
for $0<\theta<\pi.$
Taking the real and the imaginary parts of the above formula, we get
\begin{align}\label{eq:Re}
R(\theta)\cos\Theta(\theta)+\Theta'(\theta)
&=\frac c2\log\cot\frac{\theta}2+1, \\
\label{eq:Im}
R(\theta)\sin{\Theta(\theta)}-\frac{R'(\theta)}{R(\theta)}
&=\frac{c\pi}{4}.
\end{align}
Since $q_c\prec h_c,$ we note that
$\Im{q_{c}(e^{i\theta})}=R(\theta)\sin{\Theta(\theta)}\leq c\pi/4.$
Hence, $R'(\theta)\leq 0,$ which means that
$R(\theta)$ is non-increasing in $0<\theta<\pi.$
In particular, $q_c(-1)=R(\pi)\le R(\theta)\le R(0)=q_c(1).$
The same is true for $-\pi<\theta<0$ by the symmetry.
The maximum modulus principle now implies the desired inequalities.
In particular, we note that $q_c$ is bounded on $\D.$

Now we assume that $0<c\le c_1.$
Then Lemma \ref{lem:MM} implies that $q_c$ is univalent.
By Proposition \ref{prop:prop}, we see that $q_c$ meromorphically continues to $\bD\setminus\{1,-1\}.$
Since $q_c$ is bounded, there is no pole of $q_c$ on $\bD.$
Hence, $q_c$ is analytic on $\bD\setminus\{1,-1\}.$
We next show that $R(\theta)$ is strictly decreasing in $0<\theta<\pi.$
If not, since $R(\theta)$ is non-increasing,
there is an interval $I=(a,b)$ with $0<a<b<\pi$ such that
$R(\theta)$ is constant, say $R_0,$ on $I.$
Then \eqref{eq:Im} yields that $R_0\sin\Theta(\theta)=c\pi/4$
for $\theta\in I,$ which implies that $\Theta$ is also constant on $I.$
By the identity theorem, $q_c$ must be constant, which is a contradiction.
We have proved that $R(\theta)$ is strictly decreasing in $0<\theta<\pi.$
By symmetry, the same is true for $R(-\theta)$ with $0<\theta<\pi.$

We next show that $q_c'(e^{i\theta})\ne0$ for $0<\theta<\pi.$
If not, $q_c'(\zeta_0)=0$ for some $\zeta_0=e^{i\theta_0},~0<\theta_0<\pi,$
which leads to $q(\zeta_0)=h(\zeta_0)$ by \eqref{eq:q}.
Since $q_c(z)=q_c(\zeta_0)+\psi(z)^k$ near $z=\zeta_0$
for an analytic function $\psi(z)$ with $\psi(\zeta_0)=0,~\psi'(\zeta_0)\ne0$
and an integer $k\ge2,$ we can see that $q_c(\D)$ covers a small sector with
a tip at $q_c(\zeta_0)$ and opening angle is nearly $k\pi\ge 2\pi.$
This, however, contradicts the fact that $q_c(\D)$ is contained in $h_c(\D)$
whose bounadry contains a line passing through $h_c(\zeta_0)=q_c(\zeta_0).$
By symmetry, we now see that $q_c'(\zeta)\ne0$ for $\zeta\in\partial\D\setminus
\{1,-1\}.$

Since $q_c$ maps the real interval $(-1,1)$ into the positive real axis,
the upper (lower) half-disk is mapped into the upper (lower) half-plane
by $q_c.$
Therefore, we conclude that $q_c(\D)$ is a Jordan domain and
$q_c(\partial\D\setminus\{1,-1\})$ is real analytically smooth.
By Proposition \ref{prop:prop} (iii), the curve $q_c(\partial\D)$ has
continuously varying tangent at $\zeta=\pm1$ as well.
Thus we get the last conclusion.
\end{pf}

We define $c_\alpha$ for $0<\alpha\le1$ as the largest possible number so that
$$
0<c<c_\alpha\quad \Rightarrow \quad |\arg q_c|<\frac{\pi\alpha}2~\text{on}~\D.
$$
Obviously, when $\alpha=1$ this number agrees with $c_1$ defined before
Theorem \ref{thm:qc}.
The following result reduces the computation of $\gamma(\SS{\alpha})$
to the investigation of mapping properties of the function $q_c.$

\begin{lem}\label{lem:ca}
For $0<\alpha\le1,$ the relation $\gamma(\SS{\alpha})=\pi c_\alpha/4$ holds.
\end{lem}

\begin{pf}
Let $\gamma=\gamma(\SS{\alpha}).$
Then $|\Im [zf_c''/f_c']|=|\Im h_c|<\pi c/4\le\gamma$ for $c\le 4\gamma/\pi.$
By the definition of the number $\gamma(\SS{\alpha}),$ we obtain
$f_c\in\SS{\alpha},$ which means that $|\arg q_c|<\pi\alpha/2.$
Therefore, we have $c_\alpha\ge 4\gamma(\SS{\alpha})/\pi.$

Next assume that $|\Im [zf''/f']|<\pi c_\alpha/4.$
Then $\varphi_f=1+zf''/f'\prec h_{c_\alpha}.$
We note that $c_\alpha\le c_1\le c_*.$
By Theorem \ref{thm:qc} together with Lemma \ref{lem:MM},
we see that $q_f=zf'/f\prec q_{c_\alpha}.$
Since $|\arg q_{c_\alpha}|<\pi\alpha/2,$ we have $f\in\SS{\alpha}.$
Hence, $\gamma(\SS{\alpha})\ge \pi c_\alpha/4.$
We now conclude that $\gamma(\SS{\alpha})= \pi c_\alpha/4.$
\end{pf}

\section{Proof of main results}
In order to obtain upper bounds for $c_\alpha,$ we use the
Carath\'eodory-Toeplitz theorem (see, for instance,
\cite[Theorem IV.22]{Tsuji:Potential}).

\begin{lem}\label{lem:CT}
Let $p(z)=1+b_1z+b_2z^2+\cdots$ be a formal power series.
Then, $p$ represents an analytic function on $\D$ with $\Re p>0$
if and only if
\begin{equation*}
\Delta_n(p):=
\begin{vmatrix}
2& b_1& b_2 &\cdots& b_n\\
\overline{b_1} & 2 & b_1 &\cdots & b_{n-1}\\
\overline{b_2} & \overline{b_1} & 2 &\cdots & b_{n-2}\\
\vdots&\vdots&\vdots&\ddots&\vdots\\
\overline{b_n} & \overline{b_{n-1}} & \overline{b_{n-2}} &\cdots & 2
\end{vmatrix}\ge0
\end{equation*}
for all $n\ge1.$
\end{lem}

We are now ready to prove Theorem \ref{thm:starlike}.

\begin{pf}[Proof of Theorem \ref{thm:starlike}]
Expand the function $q_c(z)$ in the form
$$
q_c(z)=1+\sum_{n=1}^\infty b_nz^n=1+b_1z+b_2z^2+\cdots.
$$
Then, by comparing the series expansion of the both sides of \eqref{eq:q}
(or, alternatively, via the formula \eqref{eq:convex}, by using the
relation $q_c(z)=zf_c'(z)/f_c(z)$), we obtain
\begin{align*}
q_c(z)&=1+\frac c2z+\frac{c^2}{12}z^2+\frac{c}{12}z^3+\frac{c^2(24-c^2)}{720}z^4
+\frac{c(5c^2+72)}{2160}z^5 \\
&\quad +\frac{c^2(c^4-24c^2+522)}{30240}z^6+\frac{c(1620+189c^2-7c^4)}{90720}z^7
+\cdots.
\end{align*}

By using computer algebra, we get
\begin{align*}
&\qquad \Delta_6(q_c)=2^{-21}\cdot 3^{-14}\cdot 5^{-6}\cdot 7^{-2}\\
&\cdot\left(9 c^{12}-3168 c^{10}+117032 c^8-5676096 c^6+456371280 c^4-10334615040
   c^2+62705664000\right) \\
&\cdot\left(9 c^{12}-2328 c^{10}+209872 c^8-9890976 c^6+266580720
   c^4-3412575360 c^2+15676416000\right)
\end{align*}
and find that its minimal positive root $\rho_6$ is approximately $3.1735.$
Hence, by Lemma \ref{lem:ca}, $\gamma(\es^*)=\pi c_1/4<\pi\rho_6/4
<2.4925<2.5.$

The first inequality in Theorem \ref{thm:starlike} follows from Corollary
\ref{cor:sqrt3}.
\end{pf}

We may increase the number $n=6$ when applying Lemma \ref{lem:CT}
to obtain a better upper bound.
See the last section for such attempts.

\begin{pf}[Proof of Theorem \ref{thm:ss}]
Let $0<\alpha<1.$
For the function $q(z)=[(1+z)/(1-z)]^\alpha,$ Mocanu \cite{Mocanu86A}
considered the corresponding open door function
$h(z)=q(z)+zq'(z)/q(z).$
The authors showed in \cite[Lemma 3.3]{LS15} that the image $h(\D)$
contains a parallel strip of the form $W_\gamma$ for some $\gamma>g(\alpha),$ where
\begin{align*}
g(\alpha)&=
\frac12\left[
(1+\alpha)\sqrt{1+2\sin(\pi\alpha/2)}-\frac{1-\alpha}{\sqrt{1+2\sin(\pi\alpha/2)}}
\right] \\
&=\frac{\alpha+(1+\alpha)\sin(\pi\alpha/2)}{\sqrt{1+2\sin(\pi\alpha/2)}}.
\end{align*}
The Mocanu theorem \cite[Theorem 2]{Mocanu86A} (see also \cite{LS15}),
which is a version of the Open Door Lemma, implies that
$\gamma(\SS{\alpha})\ge \gamma>g(\alpha).$

We next show that $g(\alpha)>\sqrt3\alpha.$
Since
\begin{align*}
g''(\alpha)
&=\pi\frac{8\cos(\alpha\pi/2)\sin(\alpha\pi/2)
-(1+\alpha)\pi\sin^2(\alpha\pi/2)-2\pi\sin(\alpha\pi/2)
-(1+\alpha)\pi}%
{8(1+\sin(\alpha\pi/2))^{3/2}} \\
&\quad\, -(1-\alpha)\frac{3\pi^2\cos^2(\alpha\pi/2)}%
{8(1+2\sin(\alpha\pi/2))^{5/2}} \\
&\le 
\frac{\pi(8\cos(\alpha\pi/2)-3\pi)\sin(\alpha\pi/2)}%
{8(1+\sin(\alpha\pi/2))^{3/2}} 
-(1-\alpha)\frac{3\pi^2\cos^2(\alpha\pi/2)}%
{8(1+2\sin(\alpha\pi/2))^{5/2}},
\end{align*}
it is easy to see that $g''(\alpha)<0$ for $0<\alpha<1,$
in other words, $g(x)$ is strictly concave.
Hence, we have the inequality
$g(\alpha)>g(0)+(g(1)-g(0))\alpha=\sqrt3\alpha$ for $0<\alpha<1.$

Finally, we consider the upper estimate.
Let $c\le c_\alpha.$
Then
$$
p(z)=[q_c(z)]^{1/\alpha}=1+\frac{c}{2\alpha}z
+\frac{(3-\alpha)c^2}{24\alpha^2}z^2
+\frac{4\alpha^2c+(1-\alpha)c^3}{48a^3}z^3+\cdots
$$
has positive real part on $\D.$
By Lemma \ref{lem:CT}, $\Delta_2(p)\ge0$ is necessary for
$p\in\Cara.$
A straightforward computation yields
\begin{align*}
\Delta_2(p)&=
\frac{(9-\alpha^2)c^4-288 \alpha^2c^2+2304 \alpha^4}{288 \alpha^4} \\
&=\frac{\{(3+\alpha)c^2-48 \alpha^2\}
\{(3-\alpha)c^2-48 \alpha^2\}}{288\alpha^4}.
%\left|\begin{array}{ccc}
% 2 & \dfrac{c}{2 a} & -\dfrac{(a-3) c^2}{24 a^2} \medskip  \\
% \dfrac{c}{2 a} & 2 & \dfrac{c}{2 a} \medskip \\
% -\dfrac{(a-3) c^2}{24 a^2} & \dfrac{c}{2 a} & 2 \\
%\end{array}\right|.
\end{align*}
By solving the inequality $\Delta_2(p)\ge0,$ we obtain
$c_\alpha\le 4\sqrt3\alpha/\sqrt{3+\alpha}.$
Now Lemma \ref{lem:ca} gives the desired upper bound.
\end{pf}

\begin{pf}[Proof of Theorem \ref{thm:char}]
%Recall that $c_0>0$ is the smallest value so that $\Re F(c_0)=0.$
%If there is no such $c_0,$ we put $c_0=+\infty.$
%We already defined $c_1$ as the largest possible number such that
%$\Re q_c>0$ on $\D$ for $0<c<c_1.$
%Thus, we need to show that $c_0=c_1.$
%First we show that $c_0\ge c_1.$
%If $c_0=+\infty,$ there is nothing to show.
%Thus we assume $c_0<+\infty.$
%By definition, the boundary of $q_{c_0}(\D)$ intersects
%(or touches) the imaginary axis.
%Therefore, it is clear that $c_0\ge c_1.$
%We show that $c_0\le c_1.$
Recall that $c_1$ is the largest possible number such that
$\Re q_c>0$ on $\D$ for $0<c<c_1.$
Theorem \ref{thm:qc} tells us that $q_{c_1}$ is a bounded
univalent function on $\D$ and that the boundary of $D=q_{c_1}(\D)$
does not touch the origin.
Then the argument function $\Theta(\theta)$ defined in the
proof of Theorem \ref{thm:qc} with $c=c_1$ satisfies
that $\Theta(0)=\Theta(\pi)=0$ and $0<\Theta(\theta)\le\pi/2$
for $0<\theta<\pi.$
By maximality, $\Theta(\theta_0)=\pi/2$ for some $0<\theta_0<\pi.$
Since $\pi/2$ is the possible largest value of $\Theta(\theta),$
we have $\Theta'(\theta_0)=0.$
Then, we substitute these into \eqref{eq:Re} to get
$(c_1/2)\log\cot(\theta_0/2)+1=0,$ equivalently,
$\theta_0=2\arctan(e^{2/c_1})=\theta_{c_1}.$
Hence,
$\Re F(c_1)=\Re q_{c_1}(e^{i\theta_0})=\Re[R(\theta_0)e^{i\Theta(\theta_0)}]=0.$
For $0<c<c_1,$ we have $\Re F(c)=\Re q_c(e^{i\theta_c})>0.$
Thus we conclude that $c_1$ is the smallest positive number
such that $\Re F(c_1)=0,$ that is to say, $c_1=c_0.$
\end{pf}

\section{Numerical experiments}

\begin{figure}
\begin{center}
\includegraphics{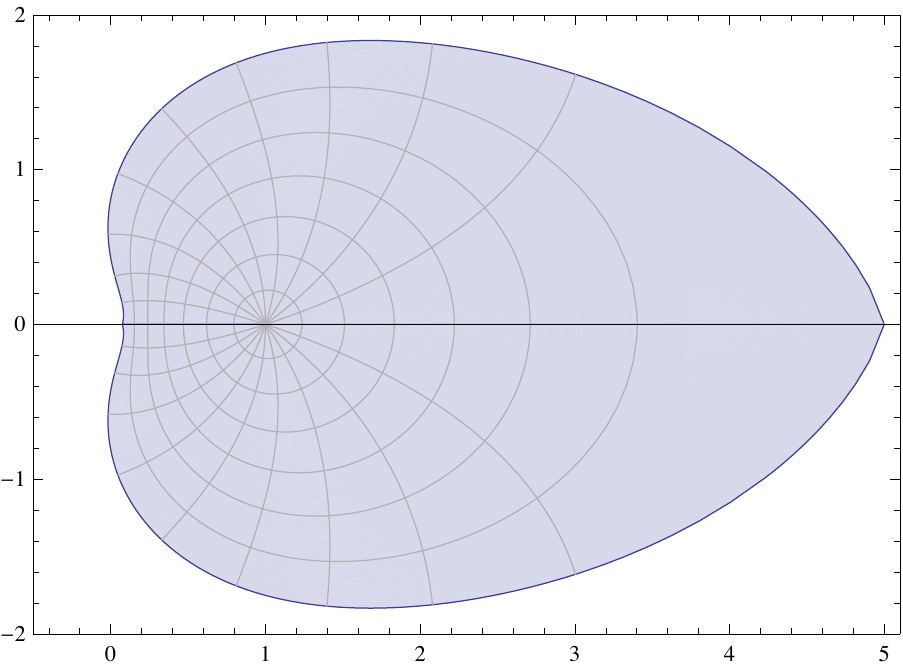}
\caption{Conformal mapping of $\D$ under $q_{c_1}.$}\label{fig:qc}
\end{center}
\end{figure}

\begin{figure}
\begin{center}
\includegraphics{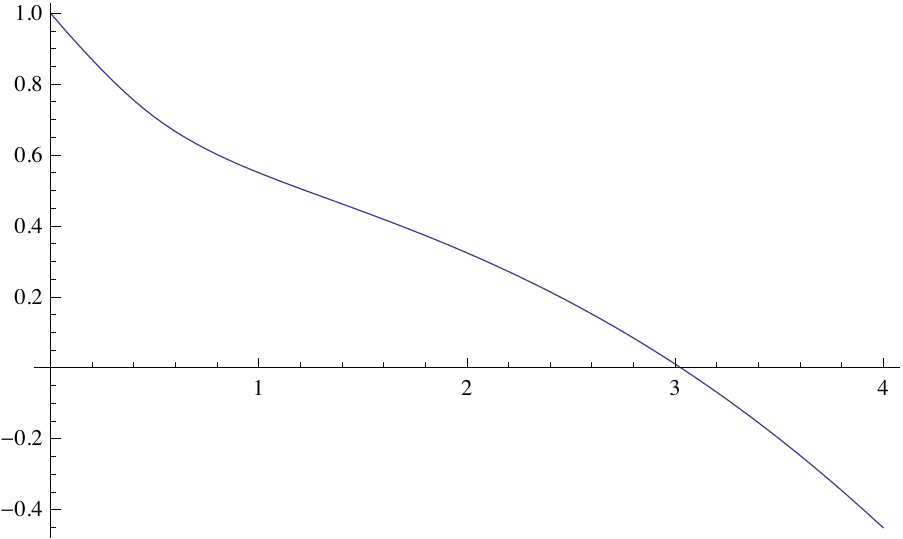}
\caption{The graph of $\Re F(c).$}\label{fig:F}
\end{center}
\end{figure}

By using Mathematica Ver.~10, we can evaluate the right-hand side of
\eqref{eq:1/q}.
In this way, we can compute the values of $q_c(z)$ numerically.
In Figure \ref{fig:F}, we exhibit the graph of the function
$c\mapsto \Re F(c),$ where $F(c)$ is given in Theorem \ref{thm:char}.
Numerical experiments give us
$\Re F(3.02756)\approx 1.06\times 10^{-6}$ and
$\Re F(3.02757)\approx -2.80\times 10^{-6}.$
Thus, if the numerical computations were correct, we would have
$c_0=c_1\approx 3.0276.$
The image of $\D$ under the mapping $q_{c_1}$ is generated in this way
(see Figure \ref{fig:qc}).

\begin{figure}
\begin{center}
\includegraphics{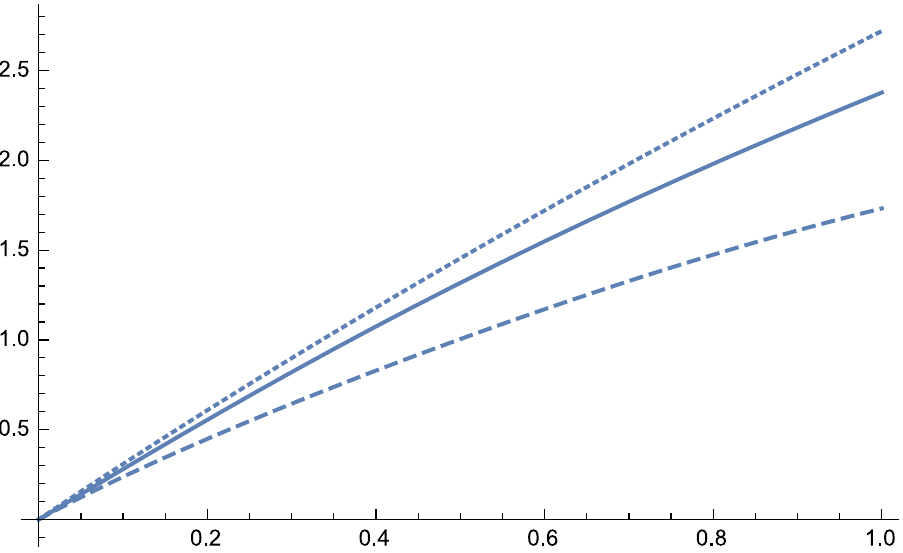}
\caption{The graphs of $\gamma(\SS{\alpha})$ (solid line) $g(\alpha)$ (dashed line) and $\sqrt{3}\pi\alpha/\sqrt{3+\alpha}$ (dotted line).}\label{fig:gamma}
\end{center}
\end{figure}

In the same way, based on Lemma \ref{lem:ca}, we can draw a graph of
the function $\alpha\mapsto\gamma(\SS{\alpha})$ together with
the upper bound $\sqrt3\pi\alpha/\sqrt{3+\alpha}$ and the lower bound
$g(\alpha)$ given in Theorem \ref{thm:ss}, see Figure \ref{fig:gamma}.
The image looks to have a corner at $q_{c_1}(\pm1).$
However, if we magnify the neighbourhood of these points large enough,
it should look smooth according to Theorem \ref{thm:qc}.

In the previous section, we obtained upper bounds for $c_1$ based on
the Carath\'eodory-Teoplitz theorem.
With the help of computer algebra, we can go further.
We can compute $\Delta_n(q_c)$ exactly as a polynomial in $c$
with rational coefficients for a small enough $n$
and find numerically the smallest positive root $\rho_n$
of the polynomial $\Delta_n(q_c)$ in $c.$
In this way, Table \ref{tb:1} can be made with the aid of Mathematica.
Thus the upper bound $2.5$ in Theorem \ref{thm:starlike} can be reduced
to some extent.
Some results are depicted in Table \ref{tb:1}.
We see that $\rho_{30}$ is close enough to the expected value
$c_1\approx 3.02756.$

\begin{center}
\begin{table}[h]%\label{tb:1}
\caption{Approximated values of $\rho_{n},$ $\gamma_{n}=\pi\rho_{n}/4$ and
$\Delta\gamma_{n}=\gamma_{n}-\gamma_{n-1}$.}\label{tb:1}
\begin{tabular}{|c|c|c|c|}
\hline
$n$ & $\rho_{n}$ & $\gamma_{n}$ & $-\Delta\gamma_{n}$\\
\hline
1&4.00000000&3.14159265&\\
2&3.46410162&2.72069905&0.42089400\\
3&3.36499696&2.64286243&0.07783660\\
4&3.33586037&2.61997861&0.02288382\\
5&3.21295295&2.52344735&0.09653126\\
6&3.17351296&2.49247125&0.03097610\\
7&3.17032183&2.48996494&0.00250631\\
8&3.13275982&2.46046381&0.02950113\\
9&3.11076636&2.44319018&0.01727363\\
10&3.10609706&2.43952292&0.00366726\\
15&3.06686241&2.40870810&0.00907899\\
20&3.04388463&2.39066140&0.00107182\\
30&3.04026630&2.38781957&0.00014363\\
\hline
\end{tabular}
\end{table}
\end{center}

\def\cprime{$'$} \def\cprime{$'$} \def\cprime{$'$}
\providecommand{\bysame}{\leavevmode\hbox to3em{\hrulefill}\thinspace}
\providecommand{\MR}{\relax\ifhmode\unskip\space\fi MR }
\providecommand{\MRhref}[2]{%
  \href{http://www.ams.org/mathscinet-getitem?mr=#1}{#2}
}
\providecommand{\href}[2]{#2}

%\bibliography{papers}
\end{document}